\documentclass{article}
\usepackage{amsmath,amssymb,amsfonts,amsthm,enumerate,xcolor,enumerate,stmaryrd,mathtools,authblk}
\usepackage[all]{xy}

\newcommand{\C}{\mathbb{C}}
\newcommand{\I}{\mathbb{I}}

\newcommand{\SE}[1]{\mathbf{SplExt}(#1)}
\newcommand{\Aut}{\textnormal{Aut}}
\newcommand{\Top}{\mathbf{Top}}
\newcommand{\Gp}{\mathbf{Gp}}
\newcommand{\TopGp}{\mathbf{TopGp}}
\newcommand{\GpAct}{\mathbf{GpAct}}
\newcommand{\Set}{\mathbf{Set}}
\newcommand{\spexc}[1]{\llbracket#1\rrbracket}
\newcommand{\op}{\textnormal{op}}
\newcommand{\SplExt}{\textnormal{SplExt}}

\newcommand{\ev}{\textnormal{ev}}
\newcommand{\two}{\mathbf{2}}

\newtheorem{theorem}{Theorem}[section]
\newtheorem{corollary}[theorem]{Corollary}
\newtheorem{proposition}[theorem]{Proposition}

\begin{document}
\title{Normalizers in the category of topological groups}
\author[1]{M.M. Clementino\thanks{The first author acknowledges partial financial support by CMUC (doi.org/10.54499/UID/00324/2025), funded by the Portuguese
Foundation for Science and Technology (FCT), grants UID/00324/2025 and UID/PRR/00324/2025.}}
\author[2]{J.R.A. Gray}
\affil[1]{CMUC, Department of Mathematics, University of Coimbra, 3000-143 Coimbra, Portugal}
\affil[2]{Mathematics Division, Department of Mathematical Sciences, Stellenbosch University, Private Bag X1, Matieland 7602, South Africa}
\maketitle
\begin{abstract}
We show that the category of morphisms of topological groups admits
split extension classifiers and hence the category of topological groups admits normalizers. 
We give explicit descriptions of both normalizers and of split extension classifiers of morphisms
of topological groups in terms of split extension classifiers of topological groups.
\end{abstract}
\section{Introduction}
It is well known that normalizers form a fundamental concept in group theory.
One of the main aims of this paper is to show that they exist in a certain form for
topological groups. Let us recall the necessary background so that we can explain
precisely what we mean. Note that for a monomorphism  $f:X\to Z$ of topological groups,
that is, an injective continuous homomorphism, if
$f(X)$ is a normal subgroup of $Z$ (as groups), then there is
unique map $\bar \gamma : X\times Z\to X$ such that for each $x\in X$ and $z\in Z$,
$f(\bar \gamma(x,z))=zf(x)z^{-1}$. In our setting, Proposition 3.1 of \cite{BOURN:2006} 
implies that $f$ is \emph{Bourn-normal} (defined below) if and only if $f(X)$ is a normal subgroup of $Z$ (as groups)
and the map $\bar \gamma$ is continuous.
We show that for an arbitrary monomorphism
$f:X\to Z$ of topological groups the set
\begin{align*}
N=\{z \in Z\mid& \text{ the auto-homeomorphism $\gamma_Z(z)$}\\
&\text{defined by $z'\mapsto zz'z^{-1}$
restricts along $f$ to $X$}\}
\end{align*}
can be equipped with a topological group structure, such that the inclusion
map $m:N\to Z$ is a monomorphism of topological groups, and the induced map $n:X\to N$
defined by $n(x)=f(x)$ (which exists because $\gamma_Z(f(x))$ restricts to $\gamma_X(x)$)
is a Bourn-normal monomorphism. Moreover, this factorization $(N,n,m)$ is the \emph{normalizer}
of $f$: it is terminal in the category of such factorizations where the requirement that the second
morphism is a monomorphism is omitted. Our result generalizes Proposition 2.13 of
\cite{BOURN_GRAY:2015} which shows that normalizers exist for those 
monomorphisms of topological groups which satisfy a certain technical condition.

It is well known that a split epimorphism $p:A\to B$ of groups is, up to isomorphism, a
semi-direct product projection, which in turn is determined by an \emph{action}
of $B$ on the kernel $X$ of $p$. It is also well known that such an action is
uniquely determined by a homomorphism $B\to \Aut(X)$. The first author together
with F.\ Cagliari \cite{CAGLIARI_CLEMENTINO:2019} showed 
(generalizing results from \cite{BORCEUX_CLEMENTINO_MONTOLI:2014}) that
essentially the same is true for topological groups, provided $\Aut(X)$ is
replaced by the group of auto-homeomorphisms of $X$ with a suitable topology.
We show here that the above facts (suitably interpreted) are true for the category of morphisms
of topological groups. We then use this result via general categorical facts
(see Proposition \ref{proposition:abstract_normalizer}, Theorem \ref{theorem:normalizer}
and Corollay \ref{corollary:normalizer} below) to construct normalizers.

\section{Preliminaries}
Recall that a morphism $n$ in a pointed category $\C$ is a normal monomorphism
if it is the kernel of some morphism. When $\C$ is moreover finitely complete
this is equivalent to $n$ being the zero class of an effective equivalence
relation; indeed if $R$ is the kernel pair of the morphism $n$ is normal
to, then $n$ forms part of a pullback
\[
\xymatrix{
N\ar[d]_{n}\ar[r]^{u} &  R\ar[d]^{\langle r_1,r_2\rangle}\\
A \ar[r]_-{\langle 1,0\rangle} & A\times A.
}
\]
We call the zero class of an arbitrary equivalence relation a Bourn-normal
monomorphism (in agreement with the definition \cite{BOURN:2000a} in the pointed context). One of our main results is to show that topological groups admit normalizers,
in the following sense: each monomorphism $f:A\to B$ admits a terminal factorization
$(N,n,u)$ where $n:A\to N$ is a \emph{Bourn-normal monomorphism} and $u:N\to B$ an arbitrary
morphism. This definition is equivalent
to the existence of normalizers in the sense of \cite{BOURN_GRAY:2015} (although we
have dropped certain additional structure in the formulation since it automatically
exists and is unique). It is
different to that appearing in \cite{GRAY:2013b} where $n$ was required to be
normal and $u$ a monomorphism.

Recall that for a pointed category $\C$, a \emph{split extension (of $B$ with kernel $X$)} is a diagram in $\C$
\begin{equation}
\label{equation: a split extension}
\xymatrix{X \ar[r]^{k} & A \ar@<0.5ex>[r]^{p} & B\ar@<0.5ex>[l]^{s}}
\end{equation}
where $k$ is the kernel of $p$, and $ps = 1_B$.  A morphism of split extensions is a diagram in $\C$

\begin{equation}
\label{equation: morphism of split ext}
\vcenter{
\xymatrix{
X \ar[r]^{k} \ar[d]_{u}& A \ar@<0.5ex>[r]^{p}\ar[d]^{v} & B\ar@<0.5ex>[l]^{s}\ar[d]^{w}\\
X' \ar[r]^{k'} & A' \ar@<0.5ex>[r]^{p'} & B'\ar@<0.5ex>[l]^{s'}
}
}
\end{equation}
where the top and bottom rows are split extensions (the domain and codomain respectively), and $vk = k' u$, $vs = s' w$ and $wp=p'v$.
Writing $\SE{\C}$ for the category of split extensions we can construct the span
\begin{equation}
\label{split ext span}
\xymatrix{
\C & \SE{\C} \ar[l]_-{P} \ar[r]^-{K} & \C
}
\end{equation}
where $P$ and $K$ are the functors sending the split extension \eqref{equation: a split extension} to $B$ and $X$, respectively. A \emph{generic split extension with kernel $X$} (in the sense
of \cite{BORCEUX_JANELIDZE_KELLY:2005a}) is a terminal object in the category $K^{-1}(X)$. We denote such a split extension as follows
\begin{equation}
\label{equation: generic split extension}
	\xymatrix{X \ar[r]^-{k_X} & H(X) \ar@<0.5ex>[r]^-{p_X} & \spexc{X}. \ar@<0.5ex>[l]^-{i_X}}
\end{equation}

For a category $\C$ with finite limits let us write $\GpAct(\C)$
for the category with objects $(B,X,\alpha)$ where $B$ and $X$ are group
objects in $\C$ and $\alpha\colon B\times X \to X$ is a morphism in $\C$
such that the diagrams
\[
\xymatrix{
B\times B\times X
\ar[r]^-{1\times \alpha}
\ar[d]_{m\times 1}
&
B\times X
\ar[d]^{\alpha}
\\
B\times X
\ar[r]_{\alpha}
&
X
}
\xymatrix{
B\times X\times X
\ar[rr]^-{\langle \pi_1, \pi_2,\pi_1,\pi_3\rangle}
\ar[d]_{\alpha\times 1}
&&
B\times X \times B\times X
\ar[d]^{\alpha\times \alpha}
\\
B\times X
\ar[r]_{\alpha}
&
X
&
X\times X
\ar[l]^{m}
}
\]
\[
\xymatrix{
B
\ar[r]^-{\langle 1,e!\rangle}
\ar[rd]_{e!}
&
B\times X
\ar[d]^{\alpha}
\\
&
X
}
\]
commute. A morphism $(f,g):(B,X,\alpha)\to (B',X',\alpha')$ is a pair of internal
group morphisms $f:B\to B'$ and $g:X\to X'$ such that the diagram
\[
\xymatrix{
B\times X
\ar[r]^-{f\times g}
\ar[d]_{\alpha}
&
B'\times X'
\ar[d]^{\alpha'}
\\
X
\ar[r]_{g}
&
X'
}
\]
commutes. Let us write $P',K':\GpAct(\C)\to \Gp(\C)$ for the functors sending a triple $(B,X,\alpha)$ to $B$ and $X$, respectively,
and defined on morphisms analogously.
\begin{theorem}\label{theorem: act}
The functor $U\colon \SE{\C}\to \GpAct(\C)$, with object map sending a split extension of the form \eqref{equation: a split extension} to $(B,X,\alpha)$ where $\alpha$ is the unique morphism
making the diagram
\[
\xymatrix@C=12ex{
B\times X
\ar[d]_{\alpha}
\ar[r]^{}^-{\langle s\pi_1,k\pi_2,s\pi_1\rangle}
&
A\times A \times A
\ar[d]
\ar[d]^{m(1\times m(1\times i))}
\\
X
\ar[r]^{k}
&
A
}
\]
commute, is an equivalence of categories. The adjoint $V$ of $U$ can be constructed via the usual semi-direct product construction.
The functors satisfy $K'U=K$, $P'U=P$, $KV=K'$, $PV=P'$.\qed
\end{theorem}

\section{Action representability of the category of morphisms of topological groups}
In this section we show that action representability of the category of topological groups
lifts to the category of morphisms of topological groups.

According to \cite{CAGLIARI_CLEMENTINO:2019} we know that the category of topological groups
has generic split extensions where $\spexc{X}=\Aut(X)$, with underlying group the group of
homeomorphisms of $X$ which are at the same time automorphisms of the group $X$, under composition.
Since in any cartesian closed category the action $\Aut(X)\times X \to X$ is \emph{evaluation}
it follows that it is evaluation in $\Top$.

For a morphism $f:X\to Z$ let us set \[\Aut(X,Z,f)=\{(\theta,\phi)\in \Aut(X)\times \Aut(Z)\mid f\theta = \phi f\}\]
with the subspace topology of $\Aut(X)\times \Aut(Z)$. $\Aut(X,Z,f)$ inherits the structure of a topological group since it
is clearly closed under composition. Let us write $q_1\colon \Aut(X,Z,f)\to \Aut(X)$, $q_2 \colon\Aut(X,Z,f)\to \Aut(Z)$ for the composite of the inclusion and
the corresponding projection of the product. We can define an action $\ev_{f}$ of $\Aut(X,Z,f)$ on $X$ as the composite
\[
\xymatrix{
\Aut(X,Z,f)\times X \ar[r]^-{q_1\times 1} & \Aut(X)\times X \ar[r]^-{\ev_X} & X.
}
\]
A straightforward element chase allows one to see that the diagram
\begin{equation}
\tag{*}
\vcenter{
\xymatrix{
\Aut(X,Z,f)\times X
\ar[d]_-{\ev_f}
\ar[r]^-{q_2\times f}
&
\Aut(Z)\times Z
\ar[d]^{\ev_Z}
\\
X
\ar[r]_{f}
&
Z
}
}
\end{equation}
commutes.
\begin{theorem}
The triple $(\Aut(X,Z,f),\Aut(Z),q_2)$ is the representing object for the functor $\SplExt_{\C^\two}(-,(X,Z,f)) : \C^{\two^{\op}}\to \Set$.
\end{theorem}
\begin{proof}
Using the equivalence of the previous section it is sufficient to show that
each morphism
\[
\xymatrix{
B\times X
\ar[d]_{\alpha}
\ar[r]^{g\times f}
&
D\times Z
\ar[d]^{\beta}
\\
X\ar[r]_{f} & Z
}
\]
of actions in $\GpAct(\C)$ admits a unique morphism to (*).
Let us write $\overline{\alpha}\colon B \to \Aut(X)$ and $\overline{\beta}\colon D\to \Aut(Z)$ for
the unique morphisms making the diagrams
\[
\xymatrix{
B\times X
\ar[d]_{\alpha}
\ar[r]^-{\overline{\alpha}\times 1}
&
\Aut(X)\times X
\ar[d]^{\ev_X}
&
D\times Z
\ar[d]_{\beta}
\ar[r]^-{\overline{\beta} \times 1}
&
\Aut(Z)\times Z
\ar[d]^{\ev_Z}
\\
X
\ar@{=}[r]
&
X
&
Z
\ar@{=}[r]
&
Z
}
\]
commute.
These morphisms have underlying maps defined by $\overline{\alpha}(b)(x)=\alpha(b,x)$
and $\overline{\beta}(d)(z)=\beta(d,z)$. Let us show that
$\langle \overline{\alpha}, \overline{\beta} \cdot g \rangle \colon B \to \Aut(X)\times \Aut(Z)$ factors through the
topological subgroup $\Aut(X,Z,f)$.
Since
\[
(\overline{\beta} \cdot g)(b)(f(x))=\beta(g(b),f(x))=f(\alpha(b,x))=(f\cdot \overline{\alpha}(b))(x)
\]
it follows that $(\overline{\beta} \cdot g)(b) \cdot f = f \cdot \overline{\alpha}(b)$ and hence
$\langle \overline{\alpha}, \overline{\beta} \cdot g \rangle$ factors through $\Aut(X,Z,f)$. Let
us write $u$ for this morphism. Note that $q_2\cdot u = \overline{\beta} \cdot g$
so that $(u,\overline{\beta}) : (B,D,g) \to (\Aut(X,Z,f),\Aut(Z),q_2)$
and so we obtain a morphism
\[
\xymatrix{
B\times X
\ar[dd]_{\alpha}
\ar[rr]^{g\times f}
\ar[rd]_{u\times 1}
&
&
D\times Z
\ar[dd]^(0.3){\beta}
\ar[rd]^{\overline{\beta} \times 1}
\\
&
\Aut(X,Z,f)\times X
\ar[dd]_(0.3){\ev_f}
\ar[rr]^(0.4){q_2\times f}
&
&
\Aut(Z)\times Z
\ar[dd]^{\ev_Z}
\\
X
\ar[rr]_(0.6){f}
\ar@{=}[rd]
&
&
Z
\ar@{=}[rd]
\\
&
X
\ar[rr]_{f}
&
&
Z
}
\]
as desired. Uniqueness follows easily from the fact that $q_1$ and $q_2$ are jointly monic.
\end{proof}

\begin{corollary}
For any small category $\I$ the category $\TopGp^\I$ is action representable.
\end{corollary}
\begin{proof}
This essentially follows from Theorem 4.2 (a) $\Rightarrow$ (b) of \cite{GRAY:2017},
which essentially depends on Proposition 2.22. Just note that
the proof and statement of Proposition 2.22 of \cite{GRAY:2017} remain
valid if ``finite'' is replaced by ``small''.
\end{proof}
\section{Normalizers}
In this section we explain how to construct normalizers in $\TopGp$.

For $\C$ a pointed category such that $\C^{\two}$ has generic split extensions,
recall from
\cite{GRAY:2017} that each generic split
extension in
$\mathbb{C}^\two$ has codomain a generic split extension in $\C$.
Let us denote the generic split extension with kernel $f:X\to Z$
in $\mathbb{C}^\two$ by
\begin{equation}\label{generic_split_ext_in_mor_cat}
\vcenter{
\xymatrix{
X\ar[d]^{f} \ar[r]^-{k_f} & 
H(X,Z,f)\ar@<0.5ex>[r]^-{p_f}\ar[d]^{\bar q_2} & 
[X,Z,f]\ar@<0.5ex>[l]^-{i_f}\ar[d]^{q_2}\\
Z \ar[r]^-{k_Z} & H(Z)\ar@<0.5ex>[r]^-{p_Z} & [Z].\ar@<0.5ex>[l]^-{i_Z}
}
}
\end{equation}
Note that in particular this means there is a unique morphism
\[
\xymatrix@!@C=-1ex@R=-2ex{
X
\ar[rr]^-{\langle 1,0\rangle}
\ar[dd]_{f}
\ar@{=}[rd]
&
&
X\times X
\ar@<0.5ex>[rr]^-{\pi_2}
\ar[dd]^(0.3){f\times f}
\ar[rd]^{\bar \gamma_f}
&
&
X
\ar@<0.5ex>[ll]^-{\langle 1,1 \rangle}
\ar[dd]^(0.3){f}
\ar[rd]^{\gamma_f}
\\
&
X
\ar[dd]^(0.3){f}
\ar[rr]^(0.3){k_f}
&
&
H(X,Z,f)
\ar@<0.5ex>[rr]^(0.4){p_f}
\ar[dd]^(0.3){\bar q_2}
&
&
[X,Z,f]
\ar@<0.5ex>[ll]^(0.6){i_f}
\ar[dd]^{q_2}
\\
Z
\ar[rr]^(0.3){\langle 1,0\rangle}
\ar@{=}[rd]
&
&
Z\times Z
\ar@<0.5ex>[rr]^(0.7){\pi_2}
\ar[rd]^{\bar \gamma_Z}
&
&
Z
\ar@<0.5ex>[ll]^(0.3){\langle 1,1 \rangle}
\ar[rd]^{\gamma_Z}
\\
&
Z
\ar[rr]^-{k_Z}
&
&
H(Z)
\ar@<0.5ex>[rr]^-{p_Z}
&
&
[Z].
\ar@<0.5ex>[ll]^-{i_Z}
}
\]
\begin{proposition}\label{proposition:abstract_normalizer}
Let $\C$ be a pointed protomodular category such that $\C^{\two}$ admits generic split extensions. For a monomorphism $f:X\to Z$,
the object $(N,n,m)$ constructed from the diagram
\[
\xymatrix{
X
\ar@/^3ex/[rrd]^{\gamma_f}
\ar@/_3ex/[rdd]_{f}
\ar[rd]^{n}
\\
&
N
\ar[d]_{m}
\ar[r]^-{u}
&
\spexc{X,Z,f}
\ar[d]^{q_2}
\\
&
Z
\ar[r]_-{\gamma_Z}
&
\spexc{Z}
}
\]
in which the square is a pullback and the morphism $n$ is the unique morphism making the diagram commute,
is the normalizer of $f$.
\end{proposition}
\begin{proof}
Consider the diagram
\[
\xymatrix@!@C=-1ex@R=-2ex{
X
\ar[rr]^-{l}
\ar[dd]_{f}
\ar@{=}[rd]
&
&
T
\ar@<0.5ex>[rr]^-{r}
\ar[dd]^(0.3){\bar m}
\ar[rd]^{\bar u}
&
&
N
\ar@<0.5ex>[ll]^-{t}
\ar[dd]^(0.3){m}
\ar[rd]^{u}
\\
&
X
\ar[dd]^(0.3){f}
\ar[rr]^(0.3){k_f}
&
&
H(X,Z,f)
\ar@<0.5ex>[rr]^(0.4){p_f}
\ar[dd]^(0.3){\bar q_2}
&
&
[X,Z,f]
\ar@<0.5ex>[ll]^(0.6){i_f}
\ar[dd]^{q_2}
\\
Z
\ar[rr]^(0.3){\langle 1,0\rangle}
\ar@{=}[rd]
&
&
Z\times Z
\ar@<0.5ex>[rr]^(0.7){\pi_2}
\ar[rd]^{\bar \gamma_Z}
&
&
Z
\ar@<0.5ex>[ll]^(0.3){\langle 1,1 \rangle}
\ar[rd]^{\gamma_Z}
\\
&
Z
\ar[rr]^-{k_Z}
&
&
H(Z)
\ar@<0.5ex>[rr]^-{p_Z}
&
&
[Z]
\ar@<0.5ex>[ll]^-{i_Z}
}
\]
obtained by forming the pullback of the bottom and front faces.
Lemma 2.6 of \cite{GRAY:2017} tells us that the front morphism is $K$-precartesian,
which by Proposition 2.6 \cite{BOURN_GRAY:2015} implies that the back morphism is $K$-precartesian.
The proof of Proposition 2.4 of \cite{BOURN_GRAY:2015} (see also Remark 4.4 of \cite{GRAY:2017})
shows that $(N,n,m)$ is the normalizer of $f$.
\end{proof}
\begin{theorem}\label{theorem:normalizer}
For a monomorphism $f:X\to Z$ the subset
\[N=\{(\theta,z)\in \Aut(X)\times Z\mid \gamma_Z(z) f = f \theta\},\]
with subspace topology, is a topological subgroup of $\Aut(X)\times Z$
and $(N,n,m)$, where $n(x)=(\gamma_X(x),f(x))$ and $m(\theta,z)=z$,
is the normalizer of $f$.
\end{theorem}
\begin{proof}
By the previous proposition we know that the normalizer can be constructed from the pullback
\[
\xymatrix{
N
\ar[d]_{m}
\ar[r]^-{u}
&
\Aut(X,Z,f)
\ar[d]^{q_2}
\\
Z
\ar[r]_-{\gamma_Z}
&
\Aut(Z).
}
\]
Since the lower square and the outer arrows in the diagram
\[
\xymatrix{
N
\ar@/_10ex/[dd]_{m}
\ar[r]^-{u}
\ar[d]^{v}
&
\Aut(X,Z,f)
\ar@/^10ex/[dd]^{q_2}
\ar[d]_{\langle q_1,q_2\rangle}
\\
\Aut(X)\times Z
\ar[d]^{\pi_2}
\ar[r]_-{1\times \gamma_Z}
&
\Aut(X)\times \Aut(Z)
\ar[d]_{\pi_2}
\\
Z
\ar[r]_{\gamma_Z}
&
\Aut(Z)
}
\]
form pullbacks, there is a unique morphism $v$ making the top square into a pullback.
This completes the proof of the claim since our description of $N$ follows
from standard facts about pullbacks.
\end{proof}
\begin{corollary}\label{corollary:normalizer}
For a monomorphism $f:X\to Z$ the subset
\[N'=\{z\in Z\mid \gamma_Z(z) \text{ restricts along $f$ to an auto-homeomorphism of $X$}\}\]
can be equipped with a topological group structure such that the inclusion
$m'$ of $N'$ in $Z$ is a monomorphism and there is a induced morphism $n':X\to N'$
such that $(N',n',m')$ is the normalizer of $f$.
\end{corollary}
\begin{proof}
Since $m$ of the previous theorem is a monomorphism, there is 
a topological group isomorphism $h:N\to m(N)$. Now just note that $N'=m(N)$.
\end{proof}
\bibliography{myrefs}{}\bibliographystyle{mybibstyle}
\end{document}